# Routing Model for Multicommodity Freight in an Intermodal Network under Disruptions

M. Majbah Uddin
Department of Civil and Environmental Engineering
University of South Carolina
300 Main Street, Columbia, SC 29208
Email: muddin@cec.sc.edu

and

Nathan Huynh (*Corresponding author*)
Department of Civil and Environmental Engineering
University of South Carolina
300 Main Street, Columbia, SC 29208
Telephone: (803) 777-8947
Fax: (803) 777-0670
Email: huynhn@cec.sc.edu

**ABSTRACT**

This paper presents a mathematical model for the routing of multicommodity freight in an intermodal network under disruptions. A stochastic mixed integer program is formulated, which minimizes not only operational costs of different modes and transfer costs at terminals but also penalty costs associated with unsatisfied demands. The Sample Average Approximation algorithm is used to solve this challenging problem. The developed model is applied to two networks, a hypothetical 15-node network and an actual intermodal network in the Gulf Coast, Southeastern and mid-Atlantic regions of the U.S., to demonstrate its applicability, with explicit consideration of disruptions at links, nodes, and terminals. The model results indicate that under disruptions, goods in the study region should be shipped via road-rail intermodal due to the built-in redundancy of the freight transport network. Additionally, the routes generated by the model are found to be more robust than those typically used by freight carriers.

**INTRODUCTION**
The freight transport network is an essential component of the economy as it supports supply chains by connecting spatially-separated origins and destinations of supply and demand. As such, it needs to be robust and resilient to support and enhance economic development. Due to the increase in international trade, freight flows have increased significantly and this trend is expected to continue in the future (*1*). For example, a daily average of 54 million tons of freight moved through the U.S. transportation system in 2012. The projected freight flows will stress both public and private infrastructures as more elements of the network reach or exceed capacity, which in turn will affect network performances (*2*).

The freight transport network is vulnerable to various disruptions. A disruptive event can be a natural disaster (e.g., earthquake, flooding, tornado, and hurricane) or a man-made disaster (e.g., accident, labor strike, and terrorism). A number of such disasters have occurred recently that severely impacted the freight transport network. The earthquake that occurred in 1994 on the Hayward Fault in San Francisco, CA caused more than 1,600 road closures and damaged most of the toll bridges and major highways (*3*). The collapse of the I-35W bridge in Minneapolis affected about 140,000 daily vehicle trips and the daily re-routing cost was \$400,000 for the impacted users (*4*). The Mississippi river flooding in 2011 impacted a major freight route, I-40, in Arkansas. Hurricane Sandy made landfall over the New York and New Jersey region in 2012 caused billions of dollars in damage and severely flooded streets and tunnels along the East Coast of the U. S. Due to the labor strike at the Port of Long Beach in 2012, the movement of \$650 million worth of goods was halted each day (*5*). These events highlight that damage to the transportation network not only disrupt transportation services but also result in economic losses and sociological effects. Disruptions in freight movements have a number of ramifications: (1) receivers will not receive their goods on time, (2) carriers need to find alternative routes to transport the goods that are impeded by the disruption, and (3) shippers need to adjust their supply chains to account for the disruption. For these reasons, adequate redundancy in the freight transport network is needed to prevent significant service losses in the event of a disruption (*6*).

This paper proposes a stochastic model for the routing of multicommodity freight on a road-rail intermodal network that is subject to various disruptions. The traditional intermodal location-routing model is extended to take into account potential network disruptions. The model can be used by carriers to determine the optimal road segments (highway links), rail segments (rail lines), and intermodal terminals to use under different types of disruptions. Since the exact evaluation of the stochastic model is difficult or impossible (*7*), the developed model is solved using the Sample Average Approximation algorithm proposed by Santoso et al. (*8*).

**LITERATURE REVIEW**
The multimodal freight transportation planning problem has been studied by many researchers over the past few decades, and its study was accelerated during the last decade (*9*). One of the earlier studies was done by Crainic and Rousseau (*10*), which presented a general modeling and algorithmic framework for the multicommodity, multimode freight service network to be used at the strategic and tactical planning level. The objective of their model is to minimize costs and delays, if a single authority controls the supply of transportation services and routing of goods through the service network. Their model considered capacitated network elements (i.e., roadways, rail lines, and terminals have finite capacities) and a penalty cost for excess assignment over capacity.

The majority of the studies that deal with intermodal freight shipments seek to minimize routing cost. Barnhart and Ratliff (*11*) proposed a model for minimizing routing cost in a road-rail intermodal network. Their model was to help shippers in deciding routing options. It used shortest path and matching algorithmic procedures to achieve the objective. Boardman et al. (*12*) developed a software-based decision support system (DSS) to assist shippers in making the best selection given a combination of modes. The crux of this DSS is the calculation of least-cost paths using a k-shortest path method, while requiring the transportation costs of all modes and transfer costs between modes as input. A similar approach was used by Song and Chen (*13*) in their development of mode selection software. However, the modes considered by Song and Chen had pre-scheduled departure times. The authors concluded that the minimum cost delivery problem is equivalent to the shortest path problem if the release time at the origin and the due date at the destination are provided.

A number of studies have addressed the intermodal routing problem with time windows. Ziliaskopoulos and Wardell (*14*) proposed an algorithm for finding the optimal time-dependent intermodal path in a multimodal transportation network. Their algorithm considered mode and arc switching delays. Xiong and Wang (*15*) developed a bi-level multi-objective model and genetic algorithmic framework for the routing problem with time windows in a multimodal network. Ayar and Yaman (*16*) investigated an intermodal multicommodity routing problem where release times and due dates of commodities were pre-scheduled in a planning horizon.

All of the aforementioned studies assume that the freight transport network is always functioning and is never disrupted, which is not realistic. To account for natural or man-made disruptions, some researchers have studied the reliability, vulnerability, and resiliency of transportation networks. Snyder and Daskin (*17*) presented a reliable uncapacitated location problem considering failure of facilities in the network. Their model finds reliable facility location by taking into account the expected transportation cost after failure, in addition to the minimum operational cost. Cui et al. (*18*) extended this work to consider failures with site-dependent probabilities and re-routing of customers when there are failures. Peng et al. (*19*) also considered disruptions of facility in reliable logistics network design. Their mixed integer program not only minimizes the nominal cost but also reduces disruption risks by employing the *p*-robustness criterion.

A resilient freight transport network is one that can recover from any disruption by preventing, absorbing, or mitigating its effects (*9*). A decision model to address disruptive events in an intermodal freight transport network was proposed by Huang et al. (*20*). Their model re-routes flows if the forecasted delay on a distressed link exceeds a pre-specified threshold. In a study performed by Chen and Miller-Hooks (*21*), a method to quantify resilience of an intermodal freight transport network was developed. They formulated a stochastic mixed integer program that aims to minimize unsatisfied demands during disruptions. Their model was solved using several exact algorithms; however, the application was limited to only small-scale networks due to high computational time requirements. Miller-Hooks et al. (*22*) extended this work to maximize freight transport network resiliency by implementing preparedness and recovery activities within a given budget. A stochastic program was developed which maximizes freight flows in the network under disruptions. Similar to their previous study, the model was applied to the same small-scale networks.

A few studies have considered network vulnerability in the planning decision. Peterson and Church (*23*) investigated rail network vulnerability by formulating both uncapacitated and capacitated routing-based model, and applied their model to a statewide network. Garg and Smith (*24*) presented a methodology for designing a survivable multicommodity flow network. Their model analyzes failure scenarios involving multiple arcs. Most recently, Gedik et al. (*25*) assessed network vulnerability and re-routing of coal by rail when disruptions occur in the network.

This study fills a gap in the literature by addressing the multicommodity routing problem in an intermodal road-rail network that is subject to disruptions. This study is most closely related to the works performed by Chen and Miller-Hooks (*21*) and Miller-Hooks et al. (*22*) in that they focus on solving the road-rail intermodal freight routing problem with explicit consideration of network disruptions. However, there are several notable differences between our work and theirs: (1) our study considers the multicommodity aspect (different commodities may have different delivery requirements and some commodities might need to be separated to facilitate early or delayed delivery); (2) our study proposes a new model that uses a link-based formulation; and (3) our model is applied to an actual large-scale intermodal freight network.

## MODEL FORMULATION

The formulation assumes that a road-rail intermodal freight transportation network is represented by a directed graph $G = (N, A)$, where $N$ is the set of nodes and $A$ is the set of links. Set $N$ consists of the set of major highway intersections $H$, the set of major rail junctions $R$, and the set of intermodal terminals $S$, i.e., $N = H \cup R \cup S$. Set $A$ consists of the set of highway links $A_h$ and the set of railway links $A_r$, i.e., $A = A_h \cup A_r$. Shipments can change mode at the intermodal terminal nodes $S$. Each highway link $(i, j) \in A_h$ and railway link $(i, j) \in A_r$ have unit transportation costs associated with them for each commodity $k \in K$ shipment. Each intermodal terminal $s \in S$ has also a unit transfer cost for each commodity $k \in K$ shipment. Another important cost parameter is the penalty cost of unsatisfied demand $\Psi$. The capacity of each highway link, railway link, and intermodal terminal are disruption-scenario dependent, i.e., capacities will be different at different disruption scenarios. Similarly, the travel time on highway and railway links and the transfer time at terminals is disruption-scenario dependent.

**Sets / Indices**

| | |
|---|---|
| $H$ | set of major highway intersections |
| $R$ | set of major rail junctions |
| $S$ | set of candidate intermodal terminals |
| $A_h$ | set of highway links |
| $A_r$ | set of railway links |
| $C$ | set of OD pairs |
| $K$ | set of commodities |
| $P^c$ | set of paths $p$ connecting OD pair $c$ |
| $\Omega$ | set of disruption scenarios |
| $k$ | commodity type, $k \in K$ |
| $i, j, s$ | node, $i, j, s \in N$ |
| $c$ | an OD pair, $c \in C$ |
| $\omega$ | a disruption scenario, $\omega \in \Omega$ |

**Parameters**

| | |
|---|---|
| $d_k^c$ | original demand of commodity $k \in K$ between OD pair $c \in C$ |
| $\Psi$ | unit penalty cost for unsatisfied demand |
| $\beta_{ijk}$ | unit cost of transporting commodity $k \in K$ by truck in link $(i, j) \in A_h$ |
| $\tilde{\beta}_{ijk}$ | unit cost of transporting commodity $k \in K$ by rail in link $(i, j) \in A_r$ |
| $\beta_{sk}$ | unit cost of transferring commodity $k \in K$ in intermodal terminal $s \in S$ |
| $Q_{ij}(\omega)$ | capacity of highway link $(i, j) \in A_h$ under disruption $\omega$ |
| $\tilde{Q}_{ij}(\omega)$ | capacity of railway link $(i, j) \in A_r$ under disruption $\omega$ |
| $Q_s(\omega)$ | capacity of intermodal terminal $s \in S$ under disruption $\omega$ |
| $t_{ij}(\omega)$ | travel time on highway link $(i, j) \in A_h$ under disruption $\omega$ |
| $\tilde{t}_{ij}(\omega)$ | travel time on railway link $(i, j) \in A_r$ under disruption $\omega$ |
| $t_s(\omega)$ | processing time in intermodal terminal $s \in S$ under disruption $\omega$ |
| $T_k^c$ | delivery time for commodity $k \in K$ between OD pair $c \in C$ |
| $\Xi$ | sufficiently large number |
| $\varepsilon$ | sufficiently small number |

**Continuous Variables**

$X_{ijk}^c(\omega)$          fraction of commodity $k \in K$ transported in highway link $(i,j) \in A_h$ between OD pair $c \in C$ under disruption $\omega$

$\tilde{X}_{ijk}^c(\omega)$          fraction of commodity $k \in K$ transported in railway link $(i,j) \in A_r$ between OD pair $c \in C$ under disruption $\omega$

$U_k^c(\omega)$          unsatisfied demand of commodity $k \in K$ between OD pair $c \in C$ under disruption $\omega$

$F_{sk}^c(\omega)$          fraction of commodity $k \in K$ between OD pair $c \in C$ transferred at terminal $s \in S$ under disruption $\omega$

**Indicator Variables**

$Y_{sk}^c(\omega)$          binary variable indicating whether or not intermodal terminal $s \in S$ is selected for commodity $k \in K$ between OD pair $c \in C$ under disruption $\omega$ ($= 1$ if intermodal terminal $s$ is selected for commodity $k$ between OD pair $c$ , $= 0$ otherwise)

$\delta_{ijk}^c(\omega)$          binary variable indicating whether or not there is any flow in highway link $(i,j) \in A_h$ for commodity $k \in K$ between OD pair $c \in C$ under disruption $\omega$ ($= 1$ if highway link $(i,j)$ carries flow of commodity $k$ between OD pair $c$ , $= 0$ otherwise)

$\tilde{\delta}_{ijk}^c(\omega)$          binary variable indicating whether or not there is any flow in railway link $(i,j) \in A_r$ for commodity $k \in K$ between OD pair $c \in C$ under disruption $\omega$ ($= 1$ if railway link $(i,j)$ carries flow of commodity $k$ between OD pair $c$ , $= 0$ otherwise)

**Model Formulation**

The stochastic multicommodity intermodal freight shipment routing (**SMIFR**) problem is formulated as follows.

$$\text{Min} \quad \mathrm{E}_\omega \left[ \sum_{k \in K} \sum_{c \in C} \left( d_k^c \left( \sum_{(i,j) \in A_h} \beta_{ijk} X_{ijk}^c(\omega) + \sum_{(i,j) \in A_r} \tilde{\beta}_{ijk} \tilde{X}_{ijk}^c(\omega) + \sum_{s \in S} \beta_{sk} F_{sk}^c(\omega) \right) + \Psi U_k^c(\omega) \right) \right] \quad (1)$$

Subject to

$$\sum_{(i,m) \in A_h} X_{imk}^c(\omega) - \sum_{(m,i) \in A_h} X_{mik}^c(\omega) \begin{cases} \leq & +1 \quad \text{if } i = \mathrm{ori}^c \\ \geq & -1 \quad \text{if } i = \mathrm{des}^c, \ \forall i \in H, k \in K, c \in C, \omega \in \Omega \\ = & 0 \quad \text{otherwise} \end{cases} \quad (2)$$

$$\sum_{i \in \mathrm{ori}^c} X_{imk}^c(\omega) - \sum_{j \in \mathrm{des}^c} X_{mjk}^c(\omega) = 0, \ \forall k \in K, c \in C, \omega \in \Omega \quad (3)$$

$$X_{imk}^c(\omega) \leq \delta_{imk}^c(\omega), \ \forall (i,m) \in A_h, k \in K, c \in C, \omega \in \Omega \quad (4)$$

$$X_{mik}^c(\omega) + \delta_{imk}^c(\omega) \leq 1, \ \forall (m,i) \in A_h, k \in K, c \in C, \omega \in \Omega \quad (5)$$

$$\sum_{(i,m) \in A_h} X_{imk}^c(\omega) \sum_{(m,i) \in A_h} X_{mik}^c(\omega) = 0, \ \forall i \in ori^c, k \in K, c \in C, \omega \in \Omega \quad (6)$$

$$\sum_{(i,n)\in A_r}\widetilde{X}_{ink}^c(\omega)-\sum_{(n,i)\in A_r}\widetilde{X}_{nik}^c(\omega)=0,\ \forall i\in R,k\in K,c\in C,\omega\in\Omega \tag{7}$$

$$\sum_{(s,m)\in A_h}X_{smk}^c(\omega)-\sum_{(m,s)\in A_h}X_{msk}^c(\omega)+\sum_{(s,n)\in A_r}\widetilde{X}_{snk}^c(\omega)-\sum_{(n,s)\in A_r}\widetilde{X}_{nsk}^c(\omega)=0,$$
$$\forall s\in S,k\in K,c\in C,\omega\in\Omega \tag{8}$$

$$\left(\sum_{(s,n)\in A_r}\widetilde{X}_{snk}^c(\omega)-\sum_{(n,s)\in A_r}\widetilde{X}_{nsk}^c(\omega)\right)(1-Y_{sk}^c(\omega))=0,\ \forall s\in S,k\in K,c\in C,\omega\in\Omega \tag{9}$$

$$F_{sk}^c(\omega)=\left|\sum_{(s,m)\in A_h}X_{smk}^c(\omega)-\sum_{(m,s)\in A_h}X_{msk}^c(\omega)\right|,\ \forall s\in S,k\in K,c\in C,\omega\in\Omega \tag{10}$$

$$\varepsilon Y_{sk}^c(\omega)\le F_{sk}^c(\omega)\le Y_{sk}^c(\omega),\ s\in S,k\in K,c\in C,\omega\in\Omega \tag{11}$$

$$\sum_{(i,j)\in(A_h\cap p)}\delta_{ijk}^c(\omega)t_{ij}(\omega)+\sum_{(i,j)\in(A_r\cap p)}\widetilde{\delta}_{ijk}^c(\omega)\widetilde{t}_{ij}(\omega)+\sum_{s\in(S\cap p)}Y_{sk}^c(\omega)t_s(\omega)\le T_k^c,$$
$$\forall p\in P^c,k\in K,c\in C,\omega\in\Omega \tag{12}$$

$$\varepsilon\delta_{ijk}^c(\omega)\le X_{ijk}^c(\omega)\le\delta_{ijk}^c(\omega),\ \forall(i,j)\in A_h,k\in K,c\in C,\omega\in\Omega \tag{13}$$

$$\varepsilon\widetilde{\delta}_{ijk}^c(\omega)\le\widetilde{X}_{ijk}^c(\omega)\le\widetilde{\delta}_{ijk}^c(\omega),\ \forall(i,j)\in A_r,k\in K,c\in C,\omega\in\Omega \tag{14}$$

$$\sum_{k\in K}\sum_{c\in C}d_k^cX_{ijk}^c(\omega)\le Q_{ij}(\omega),\ \forall(i,j)\in A_h,\omega\in\Omega \tag{15}$$

$$\sum_{k\in K}\sum_{c\in C}d_k^c\widetilde{X}_{ijk}^c(\omega)\le\widetilde{Q}_{ij}(\omega),\ \forall(i,j)\in A_r,\omega\in\Omega \tag{16}$$

$$\sum_{k\in K}\sum_{c\in C}d_k^cF_{sk}^c(\omega)\le Q_s(\omega),\ \forall s\in S,\omega\in\Omega \tag{17}$$

$$d_k^c\left(1-\sum_{i\in H}X_{ijk}^c(\omega)\right)=U_k^c(\omega),\ \forall k\in K,c\in C,j=\mathrm{des}^c,\omega\in\Omega \tag{18}$$

$$0\le X_{ijk}^c(\omega)\le1,\ \forall(i,j)\in A_h,k\in K,c\in C,\omega\in\Omega \tag{19}$$

$$0\le\widetilde{X}_{ijk}^c(\omega)\le1,\ \forall(i,j)\in A_r,k\in K,c\in C,\omega\in\Omega \tag{20}$$

$$0\le F_{sk}^c(\omega)\le1,\ \forall s\in S,k\in K,c\in C,\omega\in\Omega \tag{21}$$

$$U_k^c(\omega)\in Z^+,\ \forall k\in K,c\in C,\omega\in\Omega \tag{22}$$

$$Y_{sk}^c(\omega)\in\{0,1\},\ \forall s\in S,k\in K,c\in C,\omega\in\Omega \tag{23}$$

$$\delta_{ijk}^c(\omega)\in\{0,1\},\ \forall(i,j)\in A_h,k\in K,c\in C,\omega\in\Omega \tag{24}$$

$$\widetilde{\delta}_{ijk}^c(\omega)\in\{0,1\},\ \forall(i,j)\in A_r,k\in K,c\in C,\omega\in\Omega \tag{25}$$

The objective function (1) seeks to minimize the total expected system cost across disruption scenarios. Specifically, the expected system cost includes the transportation cost on highway and railway links, the transfer cost at intermodal terminals, and the penalty cost for unsatisfied demands. Constraints (2)–(6) ensure flow conservation at highway nodes ( $H$ ). The notations $\text{ori}^c$ and $\text{des}^c$ denote the origin and destination node, respectively, of an OD pair $c \in C$. Similarly, constraint (7) ensures flow conservation at railway nodes ( $R$ ). Constraints (8) and (9) ensure flow conservation at intermodal terminals ( $S$ ); constraint (8) maintains the conservation if a terminal is selected whereas constraint (9) maintains conservation if the terminal is not selected. The decision variables $F_{sk}^c(\omega)$, $\forall s \in S, k \in K, c \in C, \omega \in \Omega$ are calculated in constraint (10). Constraints (9) and (10) are adapted from the work of Xie et al. (*26*). Constraint (11) establishes the relationship between decision variables $F_{sk}^c(\omega)$ and $Y_{sk}^c(\omega)$. Constraint (12) ensures that each commodity shipment is delivered before the delivery deadline $T_k^c$, $\forall k \in K, c \in C$. The relationship between decision variables $X_{ijk}^c(\omega)$ and $\delta_{ijk}^c(\omega)$ are expressed in constraint (13), and the relationship between decision variables $\widetilde{X}_{ijk}^c(\omega)$ and $\widetilde{\delta}_{ijk}^c(\omega)$ are expressed in constraint (14). Constraints (15)–(17) ensure that flows are less than or equal to the capacity of highway links, railway links, and intermodal terminals, respectively. Constraint (18) determines the unsatisfied demand $U_k^c(\omega)$, $\forall k \in K, c \in C, \omega \in \Omega$. Lastly, constraints (19)–(21) are the definitional constraints, constraint (22) is the integrality constraint, and constraints (23)–(25) are the binary constraints.

**Linear Formulation**
The proposed model is not linear, since it has several non-linear constraints: (6), (9), and (10). Non-linear models are generally very difficult to solve; thus, the non-linear constraints are reformulated to make the model tractable. The equivalent linear forms are

$$\sum_{(i,m) \in A_h} X_{imk}^c(\omega) \geq \Xi \sum_{(m,i) \in A_h} X_{mik}^c(\omega), \ \forall i \in \text{ori}^c, k \in K, c \in C, \omega \in \Omega \tag{26}$$

$$-\Xi Y_{sk}^c(\omega) \leq \sum_{(s,n) \in A_r} \widetilde{X}_{snk}^c(\omega) - \sum_{(n,s) \in A_r} \widetilde{X}_{nsk}^c(\omega) \leq \Xi Y_{sk}^c(\omega), \ \forall s \in S, k \in K, c \in C, \omega \in \Omega \tag{27}$$

$$-F_{sk}^c(\omega) \leq \sum_{(s,m) \in A_h} X_{smk}^c(\omega) - \sum_{(m,s) \in A_h} X_{msk}^c(\omega) \leq F_{sk}^c(\omega), \ \forall s \in S, k \in K, c \in C, \omega \in \Omega \tag{28}$$

Constraint (26) is equivalent to constraint (6), which prevents sub-tours. Adopting the approach used by Xie et al. (*26*), constraints (9) and (10) can be reformulated as constraints (27) and (28), respectively. By replacing constraints (6), (9) and (10) with constraints (26), (27), and (28), the revised model is a stochastic mixed integer linear program.

**ALGORITHMIC STRATEGY**
A key difficulty in solving a stochastic program is in evaluating the expectation of the objective function. One approach for accomplishing this is to approximate the expected objective function value through sample averaging. This study adopts the Sample Average Approximation (SAA) algorithm proposed by Santoso et al. (*8*). Without loss of generality, the objective function of the model can be rewritten as follows, where $\lambda$ represents the decision variables.

$$\min \, \text{E}_\omega \big[ \Theta(\lambda, \omega) \big] \tag{29}$$

**The SAA Algorithm**

**Step 1.** Generate $M$ independent disruption-scenario samples each of size $N$, i.e., $(\omega_j^1,...,\omega_j^N)$ for $j=1,...,M$. For each sample, solve the corresponding SAA problem.

$$\min \ \frac{1}{N}\sum_{n=1}^{N}\Theta(\lambda,\omega_j^n) \tag{30}$$

Let $f_N^j$ and $\hat{\lambda}_N^j$, $j=1,...,M$ be the corresponding optimal objective function value and an optimal solution of the model, respectively.

**Step 2.** Compute $\bar{f}_N$ and $\sigma_{\bar{f}_N}^2$ using the following equations.

$$\bar{f}_N := \frac{1}{M}\sum_{j=1}^{M}f_N^j \tag{31}$$

$$\sigma_{\bar{f}_N}^2 := \frac{1}{M(M-1)}\sum_{j=1}^{M}(f_N^j - \bar{f}_N)^2 \tag{32}$$

Here $\bar{f}_N$ provides a lower statistical bound for the optimal value $f*$ of the true problem, and $\sigma_{\bar{f}_N}^2$ is an estimate of the variance of the estimator.

**Step 3.** Choose a feasible solution $\tilde{\lambda}$ from the above computed solutions $\hat{\lambda}_N^j$, and generate another $N'$ independent disruption-scenario samples, i.e., $\omega^1,...,\omega^{N'}$. Then estimate the true objective function value $\tilde{f}_{N'}(\tilde{\lambda})$ and variance of this estimator as following.

$$\tilde{f}_{N'}(\tilde{\lambda}) := \frac{1}{N'}\sum_{n=1}^{N'}\Theta(\tilde{\lambda},\omega^n) \tag{33}$$

$$\sigma_{N'}^2(\tilde{\lambda}) := \frac{1}{N'(N'-1)}\sum_{n=1}^{N'}[\Theta(\tilde{\lambda},\omega^n)-\tilde{f}_{N'}(\tilde{\lambda})]^2 \tag{34}$$

In solving SAA problems, typically, $N'$ is much larger than the sample size $N$.

**Step 4.** Compute the optimality gap of the solution and variance of the gap estimator.

$$\text{gap}(\tilde{\lambda}) := \tilde{f}_{N'}(\tilde{\lambda}) - \bar{f}_N \tag{35}$$

$$\sigma_{\text{gap}}^2 = \sigma_{N'}^2(\tilde{\lambda}) + \sigma_{\bar{f}_N}^2 \tag{36}$$

## NUMERICAL EXPERIMENTS

To assess the applicability of the proposed model (**SMIFR**) and solution algorithm, two sets of experiments are conducted. The first set involves a hypothetical small-sized network with 15 nodes and 5 OD pairs. The second set involves an actual large-scale freight transport network, consisting of major highways, Class I railroads, and TOFC/COFC (Trailer on Flat Car/Container on Flat Car) intermodal terminals in the Gulf Coast, Southeastern and mid-Atlantic regions of the U.S.

## Network and Data Description

*Hypothetical Network*

Figure 1 shows the hypothetical 15-node road-rail freight transport network. Nodes 5, 7, 9, and 12 represent intermodal terminals, and node 8 represents a railway junction where trains can change track/route. The rest of the nodes represent highway intersections. The solid lines represent highway links, and the dashed lines represent railway links. The capacity of the links $Q_{ij}$ are assumed to have a uniform distribution (*23*), each with a specified range $[l_{ij}, u_{ij}]$ where $l_{ij}$ is the lower bound and $u_{ij}$ is the upper bound. The capacities of the intermodal terminals are also assumed to have a uniform distribution with a specified range. The demand in terms of number of shipments and delivery deadlines for each commodity between different OD pairs is provided in Table 1.

*Actual Network*

Figure 2 shows the actual road-rail freight transport network used in the second set of experiments. It is created using data provided by the Center for Transportation Analysis, Oak Ridge National Laboratory (*27*). As shown, it covers all of the states in the Gulf Coast, Southeastern and mid-Atlantic regions of the U.S.: Texas, Oklahoma, Louisiana, Alabama, Mississippi, Arkansas, Georgia, Florida, South Carolina, North Carolina, Tennessee, Kentucky, Virginia, Maryland, West Virginia, and Delaware. In all, the network has a total of 682 links (U.S. interstates and major highways and Class I railroads) and 187 nodes, including 44 intermodal terminals. Readers can refer to the work of Uddin and Huynh (*28*) for more details about the network. The Freight Analysis Zone (FAZ) centroids from the Freight Analysis Framework version 3 (FAF3) database (*29*) are treated as actual origins and destinations of commodity shipments. There are a total of 48 centroids in the study region. OD pairs are constructed from these 48 FAZ centroids, and demands are obtained from the FAF3 database. The demand data are filtered to include only those commodities typically transported via intermodal (*30*), and demands are converted into the number of TOFC/COFC containers using an average load of 40,000 lbs per container. It is assumed that all commodities need to be delivered within 7 days.

The transport cost on highways and railways are estimated to be $1.67 per mile per shipment (*31*) and $0.60 per mile per shipment (*32*), respectively. The transfer cost at intermodal terminals is estimated to be $70 per shipment (*33*). The travel times on highway and railway links are calculated using free-flow speeds. The number of potential paths between an origin and destination could be large. For that reason, after getting all the available paths between a specific OD pair, only those paths that have lengths less than or equal to five times of the corresponding minimum path length are considered in the path set. This approach is deemed reasonable because the discarded paths would not have satisfied the delivery deadline constraint.

## Disruption Types

Three types of disruptive-events are considered: (1) link disruption, (2) node disruption, and (3) intermodal terminal disruption. Link disruptions are modeled by randomly selecting several connected links and reducing their capacities by 50%. The travel times on the affected links are increased as a result of reduced capacities. Node disruptions are modeled by reducing the capacities of all links connected to the nodes by 80%. And, terminal disruptions are modeled by randomly selecting a number of terminals and reducing their capacities by 80%; thus, the transfer times at the impacted terminals will increase. It should be noted that affected links, nodes, or terminals are selected based on their vulnerability, and the severity of the disruption can be captured by the amount of capacity that is reduced. Recurring disruptions are not considered in the numerical experiments. However, these types of disruptions that occur continually over time and involving different links can easily be modeled given the generality of the model formulation and solution algorithm.

## Experimental Results

The proposed solution methodology is implemented in Python using Jetbrains PyCharm 4.0.5, and the IBM ILOG CPLEX 12.6 solver is used to solve the mixed integer program. Experiments are run on a personal computer with Intel Core i7 3.20 GHz processor and 8.0 GB of RAM.

### Hypothetical Network

To apply the SAA algorithm, the number of independent disruption-scenario samples ($M$) is set to 100, the sample size ($N$) is set to 1, and the number of large-size samples ($N'$) is set to 1,000 for all three types of disruption. With these values, the SAA method will produce a number of candidate routes per commodity per OD pair but no more than 100 ($M = 100$). Among these candidate routes, the optimal route is the one that yields the lowest optimality gap when each candidate route is applied to the 1,000 test scenarios ($N' = 1,000$).

Table 2(a) summarizes the input parameters and associated SAA results for the hypothetical network. The term "gap" denotes the optimality gap as defined in Equation 35, and $\sigma_{gap}$ denotes the standard deviation of the gap estimates as defined in Equation 36. In the case of link disruption, the average objective function value is \$92,439.62, with an optimality gap of \$540.87 and estimator standard deviation of \$17.69. The associated computation time is 17.5 minutes. Similar information is presented for the node and terminal disruption cases. Among the three types of disruption, the node disruption case results in the highest objective function value, which indicates that it has the most negative impact on freight logistics. Conversely, the terminal disruption case has the least impact. This result is counterintuitive because one would expect the terminal disruption to have the highest impact since it serves as a hub in the freight transport network. This is due to the network structure which allows commodities to be shipped via road more efficiently and less costly. In other words, terminals handle only a small percentage of the shipments, and thus, their disruptions have minimal impact on the freight logistics.

The corresponding optimal routes are presented in Table 2(b). Optimal routes are shown as a series of nodes in the direction of origin to destination. For example, the optimal route to ship commodity #1 between OD pair (1→15) in the event of link disruptions is: 1–3–5–8–12–13–15. Note that if a particular route does not have sufficient capacity to handle a particular shipment, then the remaining shipment is shipped via a second-best route. This is the case with commodity #3 between OD pair (1→11). There are two optimal routes: 1–2–6–14–13–11 (5% use this route) and 1–4–10–11 (95% use this route). It should be noted that the model places no restriction on the number of potential routes between each OD. Thus, a shipment could have several routes if there is insufficient capacity on the least-cost routes.

It is observed that since the network has very few rail links, most of the shipments are shipped via highway links. This finding corresponds to actual freight flows where the majority of freights are shipped via road. Furthermore, when highway links are disrupted, then railway links and terminals are more likely to be used. Again, this is a logical and expected result. An interesting result that highlights the usefulness of the model can be seen in the case of a node disruption for commodity #4 between OD pair (2→13). There is one optimal route, but it only contains 65% of the shipment which means that the remaining 35% failed to reach its destination (i.e., unsatisfied demand). There are no unsatisfied demands under link and terminal disruption cases.

### Actual Network

To understand the impact of disruptions on an actual road-rail intermodal network, several instances of each disruption type are considered. For link disruptions, four different instances are solved to investigate how the objective function value and computational time change with respect to the severity of the link disruption. The severity of the link disruption is modeled by the number of impacted links, which was set to 30, 60, 100, and 200 for the four instances. The results for link disruption are summarized in the first part of Table 3. The results indicate that increasing the number of OD pairs and commodities ($|K|$) will increase computational efforts. Furthermore, for a particular number of OD pairs, the objective function value increases with the number of impacted links. The computational time is unaffected by the severity level.

For node disruptions, the four instances considered have 5, 10, 20, and 40 nodes disrupted. As shown in the second part of Table 3, the objective function value and computational time increase with higher number of OD pairs and commodities. Unlike link disruption, the computational time is affected by the number of disrupted nodes. Specifically, there is a significant increase from 20 to 40 nodes for the 10 OD pairs case (154.9 minutes to 720.2 minutes).

For terminal disruptions, three instances are considered with 15, 30, and 44 terminals disrupted. The objective function value and computational time exhibit a similar trend with respect to disruption severity as the link and node disruption cases. Similar to the node disruption case, the computational time is affected by the number of disrupted terminals.

Collectively, the numerical results indicate that, under link and node disruptions, the majority of the commodity shipments are shipped via road-rail intermodal due to lower rail cost and due to the robust freight transport network. A similar finding is reported in a study done by Ishfaq (*34*) who concluded that the layout of the U.S. road-rail intermodal network and location of intermodal terminals provide sufficient redundancies to handle disruptions. When intermodal terminals are disrupted, the model indicates that commodities will be shipped via road directly. This result is expected since highway network is redundant and robust, as well as cost-effective.

In the aforementioned experiments, the unit penalty cost is assumed to be $10,000. This value is chosen to be high to ensure that the unsatisfied demand is minimized. To test the sensitivity of this parameter, experiments are performed where the unit penalty cost is set to $2,500, $5,000, and $7,500. It is found that these three values for penalty cost resulted in same amount of unsatisfied demand. The computation time is observed to increase as the penalty cost value decreases. It can be concluded that the solution is not sensitive to the unit penalty cost parameter, given that it is set to a sufficiently large value. To test the sensitivity of the delivery deadline parameter, an experiment is performed where the delivery deadline is set to 14 days. The solution, including objective function value and computation time, is found to be the same when the delivery time is 7 days. This result suggests that the majority of the shipments require less than 7 days to reach their destinations, and thus, extending the delivery deadline has no effect on the solution.

Figure 3(a) illustrates how the optimal route generated by the model for a particular commodity going from Greensboro, NC to Dallas, TX under node disruptions compares with an actual route that a carrier would use. The left part of Figure 3(a) shows the optimal route generated by the model (shown in red), and the right part shows the route that a freight carrier would use (*35*). By inspection, it is clear that the two routes are very similar to each other. However, there is one notable difference, and that is the model indicates road-rail intermodal to be optimal whereas the freight carrier chooses road-only. This discrepancy can be attributed to the fact that the carrier does not consider the potential node disruptions in the network.

Figure 3(b) illustrates how the optimal route generated by the model for a particular commodity going from Miami, FL to Houston, TX under link disruptions. By inspection, it is clear that the carrier chooses the route based on minimum travel time. The model, on the other hand, recognizes the potential link disruptions in the network and thereby chooses an intermodal route that avoids using the U.S. interstates (I-10 and I-12) through Louisiana. This is because historically this area is vulnerable to hurricanes, such as Rita and Katrina. This result illustrates the importance of considering network disruptions when selecting a route for multicommodity freight in an intermodal network.

## CONCLUSION

This paper developed a new stochastic mixed integer programming model (**SMIFR**) to determine the optimal routes for delivering multicommodity freight in an intermodal freight network that is subject to disruptions (e.g., link, node, and terminal disruptions). To solve this model, the Sample Average Approximation (SAA) algorithm is adopted. The model and solution algorithm was tested on a hypothetical 15-node network and an actual intermodal network in the Gulf Coast, Southeastern and mid-Atlantic regions of the U.S.

The numerical experiments indicated that the model is capable of finding the optimal solutions for both small and large networks. The model runtime for a hypothetical 15-node network was reasonable (less than 3 hours for all instances). Naturally, the model runtime will increase as the network gets larger, as well as for the number of OD pairs and commodities. While the computational time was affected by the severity level of node and terminal disruptions, it was unaffected by link disruption severity. The model results indicated that under disruptions, goods in the study region should be shipped via road-rail intermodal due to lower rail cost and due to the built-in redundancy of the freight transport network. Furthermore, the model indicated that for a particular number of OD pairs, the total system cost will increase as the number of disrupted elements increases. The routes generated by the model are shown to be more robust than those typically used by freight carriers because they are often selected without consideration of potential network disruptions.

**LIST OF TABLES**

**LIST OF FIGURES**

**TABLE 1  Number of Shipments and Delivery Deadlines**

| OD Pair | Commodity Index | Number of Shipments | Delivery Deadline (hours) |
|---|---|---|---|
| 1→15 | 1 | 40 | 84 |
|  | 2 | 35 | 72 |
|  | 3 | 22 | 60 |
|  | 4 | 20 | 72 |
| 1→11 | 1 | 30 | 72 |
|  | 2 | 35 | 72 |
|  | 3 | 40 | 48 |
| 2→13 | 1 | 42 | 60 |
|  | 2 | 30 | 48 |
|  | 3 | 50 | 60 |
|  | 4 | 55 | 48 |
| 15→4 | 1 | 35 | 72 |
|  | 2 | 45 | 60 |
|  | 3 | 50 | 72 |
|  | 4 | 30 | 60 |
| 14→3 | 1 | 45 | 48 |
|  | 2 | 30 | 60 |

**TABLE 2(a)  Experimental Results for Hypothetical Network**

|  | Link Disruption | Node Disruption | Terminal Disruption |
|---|---|---|---|
| $M$ | 100 | 100 | 100 |
| $N$ | 1,000 | 1,000 | 1,000 |
| CPU Time (min) | 17.5 | 178.4 | 0.9 |
| Objective Function Value (avg) | $92,439.62 | $93,152.09 | $59,419.30 |
| gap | $540.87 | $390.09 | $4.76 |
| $\sigma_{\text{gap}}$ | $17.69 | $3.80 | $0.37 |

**TABLE 2(b)  Optimal Routes for Hypothetical Network**

| OD Pair | Commodity Index | Optimal Routes | | |
|---|---|---|---|---|
|  |  | Link Disruption | Node Disruption | Terminal Disruption |
| 1→15 | 1 | 1–3–5–8–12–13–15 (100%) | 1–4–10–11–15 (100%) | 1–3–5–8–12–13–15 (100%) |
|  | 2 | 1–3–2–6–14–15 (100%) | 1–4–10–11–15 (100%) | 1–4–10–11–15 (100%) |
|  | 3 | 1–3–2–6–14–13–15 (100%) | 1–4–10–11–15 (100%) | 1–4–10–11–15 (100%) |
|  | 4 | 1–3–2–6–14–15 (100%) | 1–4–10–11–15 (100%) | 1–4–10–11–15 (100%) |
| 1→11 | 1 | 1–3–5–9–10–11 (100%) | 1–4–10–11 (100%) | 1–4–10–11 (100%) |
|  | 2 | 1–3–5–9–10–11 (100%) | 1–4–10–11 (100%) | 1–4–10–11 (100%) |
|  | 3 | 1–2–6–14–13–11 (5%) | 1–4–10–11 (100%) | 1–4–10–11 (100%) |
|  |  | 1–4–10–11 (95%) |  |  |
| 2→13 | 1 | 2–6–7–12–13 (98%) | 2–1–4–10–11–13 (33%) | 2–6–14–13 (100%) |
|  |  | 2–6–14–13 (2%) | 2–6–14–13 (67%) |  |
|  | 2 | 2–3–5–8–12–13 (100%) | 2–3–5–8–12–13 (73%) | 2–6–14–13 (100%) |
|  | 3 | 2–6–7–12–13 (100%) | 2–3–5–8–12–13 (100%) | 2–6–7–12–13 (100%) |
|  | 4 | 2–3–5–8–12–13 (84%) | 2–1–4–10–11–13 (65%) | 2–6–14–13 (100%) |
|  |  | 2–6–14–13 (16%) |  |  |
| 15→4 | 1 | 15–13–12–8–5–3–4 (100%) | 15–11–10–4 (100%) | 15–11–10–4 (100%) |
|  | 2 | 15–11–10–4 (100%) | 15–11–10–4 (100%) | 15–11–10–4 (100%) |
|  | 3 | 15–13–12–8–5–3–4 (100%) | 15–11–10–4 (100%) | 15–11–10–4 (100%) |
|  | 4 | 15–13–12–8–5–3–4 (100%) | 15–11–10–4 (100%) | 15–11–10–4 (100%) |
| 14→3 | 1 | 14–6–7–5–3 (100%) | 14–6–7–5–3 (100%) | 14–6–2–3 (100%) |
|  | 2 | 14–6–2–3 (100%) | 14–6–7–5–3 (100%) | 14–6–2–3 (100%) |

**TABLE 3  Experimental Results for Actual Network**

| OD | $\|K\|$ | Link Disruption | | | Node Disruption | | | Terminal Disruption | | |
|---|---|---|---|---|---|---|---|---|---|---|
| | | Impacted Link # | Obj. Func. ($, thousands) | CPU (min) | Impacted Node # | Obj. Func. ($, thousands) | CPU (min) | Impacted Terminal # | Obj. Func. ($, thousands) | CPU (min) |
| 5 | 9 | 30 | 556.4 | 29.8 | 5 | 557.7 | 29.9 | 15 | 550.2 | 29.4 |
| | | 60 | 562.6 | 29.9 | 10 | 566.0 | 29.3 | 30 | 560.6 | 29.5 |
| | | 100 | 573.9 | 29.9 | 20 | 604.1 | 33.6 | 44 | 651.9 | 36.5 |
| | | 200 | 650.5 | 29.9 | 40 | 640.9 | 115.1 | | | |
| 10 | 21 | 30 | 959.8 | 146.8 | 5 | 942.2 | 142.7 | 15 | 930.3 | 147.0 |
| | | 60 | 965.9 | 146.8 | 10 | 966.9 | 141.7 | 30 | 959.7 | 143.5 |
| | | 100 | 979.3 | 142.6 | 20 | 1,015.1 | 154.9 | 44 | 1,077.8 | 173.2 |
| | | 200 | 1,085.1 | 146.2 | 40 | 1,061.4 | 720.2 | | | |
| 20 | 43 | 30 | 1,478.8 | 484.1 | 5 | 1,461.2 | 486.2 | 15 | 1,481.3 | 483.5 |
| | | 60 | 1,484.9 | 487.8 | 10 | 1,505.8 | 479.6 | 30 | 1,534.3 | 493.7 |
| | | 100 | 1,500.9 | 486.7 | 20 | 1,558.4 | 528.2 | 44 | 1,705.8 | 636.5 |
| | | 200 | 1,625.2 | 485.5 | 40 | 1,609.5 | 1,204.9 | | | |
| 50 | 87 | 30 | 3,885.8 | 1,937.8 | 5 | 3,870.7 | 1,937.2 | 15 | 3,959.0 | 1,953.7 |
| | | 60 | 3,895.9 | 1,930.4 | 10 | 3,983.3 | 1,945.6 | 30 | 4,137.3 | 2,036.3 |
| | | 100 | 3,952.3 | 1,926.4 | 20 | 4,062.6 | 2,007.6 | 44 | * | * |
| | | 200 | 4,173.1 | 1,956.0 | 40 | * | * | | | |

*Program terminated due to memory limitation

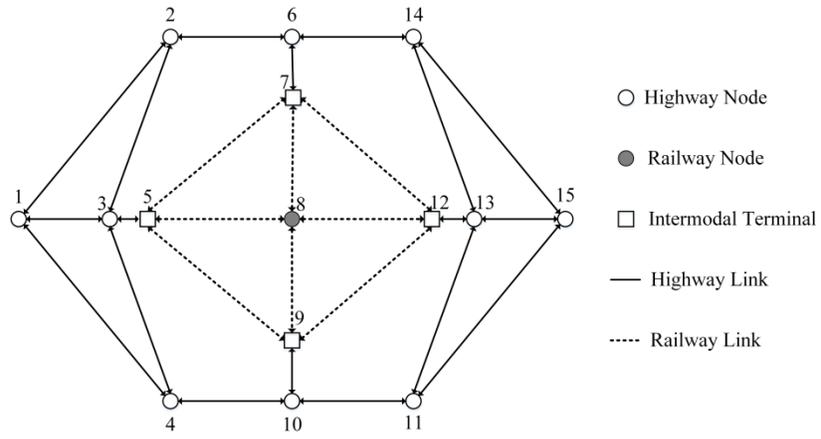

**FIGURE 1  A hypothetical 15-node road-rail freight transport network.**

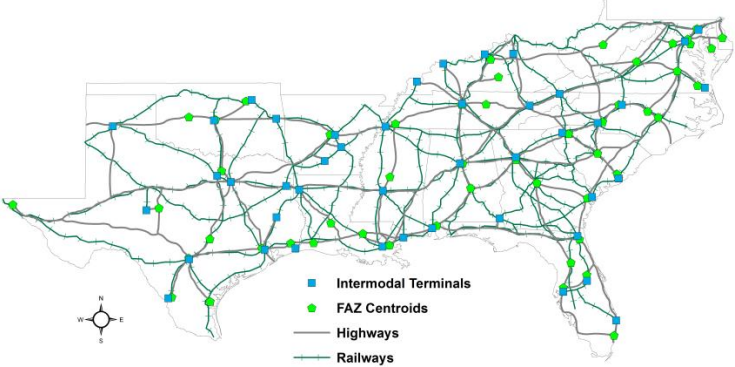

**FIGURE 2  Large-scale U. S. road-rail intermodal network.**

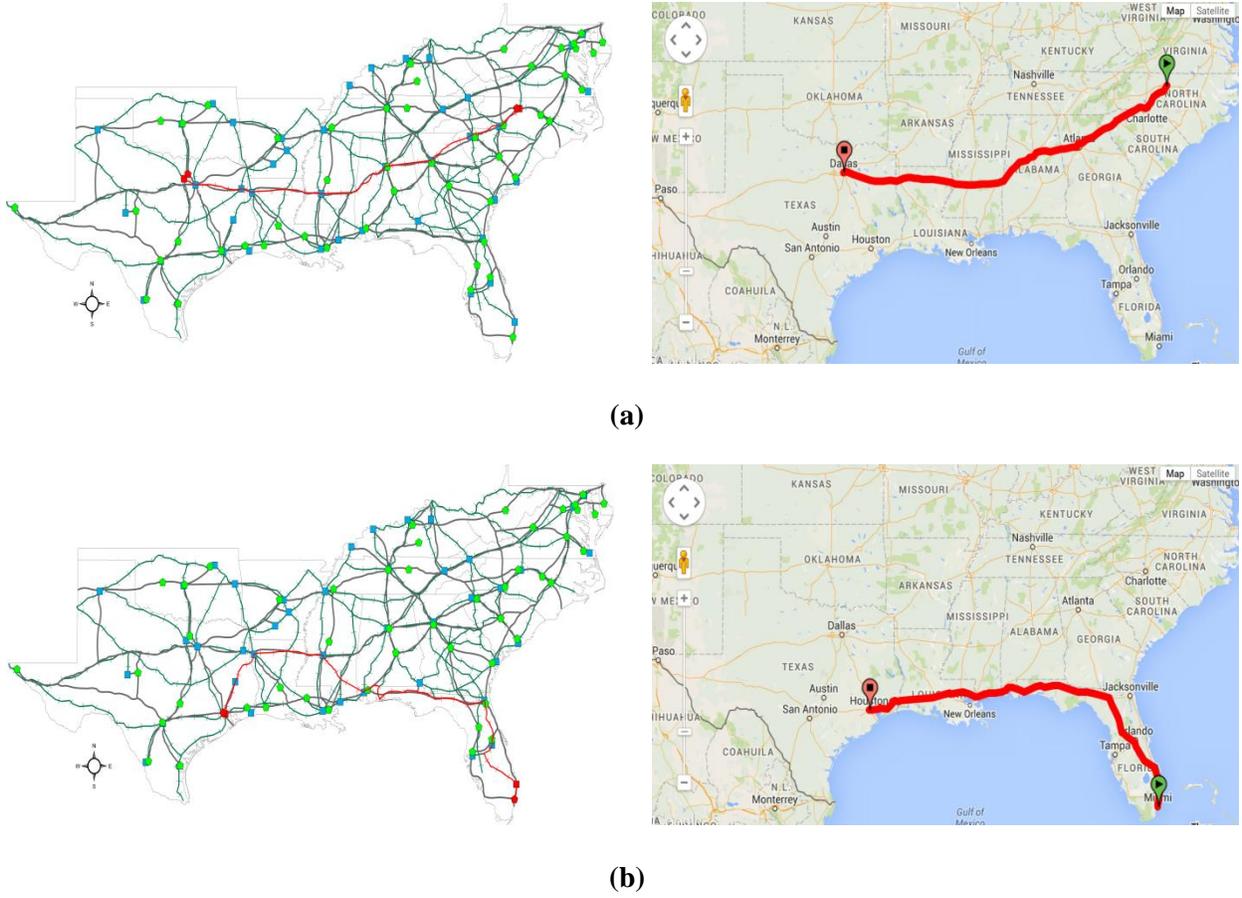

**(a)**

**(b)**

**FIGURE 3  Optimal routes for selected OD pairs: (a) Greensboro–Dallas and (b) Miami–Houston.**